\documentclass[letterpaper, 10pt, conference]{ieeeconf}
\IEEEoverridecommandlockouts                              
\usepackage[utf8]{inputenc}
\usepackage{epsfig}
\usepackage{amsmath}

\usepackage{subfigure}
\usepackage{amsfonts}
\usepackage{amsmath}
\usepackage{eso-pic}
\usepackage{textcomp}
\usepackage{comment}
\usepackage{graphicx}
\usepackage{amssymb}
\usepackage{gensymb}
\usepackage[font=small,skip=0pt]{caption}
\let\labelindent\relax
\usepackage[inline]{enumitem}
\usepackage[acronym]{glossaries} 
\usepackage{authblk}
\makeglossaries
\graphicspath{{Figures/}}
\usepackage{draftwatermark}

\makeatletter
\let\NAT@parse\undefined
\makeatother
\usepackage{hyperref}
\hypersetup{
colorlinks=true, 
breaklinks=true, 
urlcolor= black, 
citecolor=black,
linkcolor= blue, 
bookmarksopen=true,
pdftitle={Autonomous driving using an optimized neural network based adaptive LPV-MPC controller}, 
pdfauthor={Y. Kebbati},
}

\title{\LARGE \bf Optimized self-adaptive PID speed control for autonomous vehicles
}

\author[1,*]{Yassine Kebbati}
\author[1]{Naima Ait-Oufroukh}
\author[1,2]{Vincent Vigneron}
\author[1]{Dalil Ichalal}
\author[3]{Dominique Gruyer}
\affil[1]{IBISC-EA4526, Université Evry, Université Paris-Saclay, France} 
\affil[2]{School of Applied Sciences (FCA),UNICAMP, Limeira, Brazil} 
\affil[3]{IFSTTAR/ PICS-L, Université Gustave Eiffel, France \authorcr Email$^{*}$: {\tt \href{mailto:yassine.kebbati@univ-evry.fr} yassine.kebbati@univ-evry.fr}\vspace{1.5ex}}

\newcommand\AtPageUpperMyright[1]{\AtPageUpperLeft{%
 \put(\LenToUnit{0.08\paperwidth},\LenToUnit{-1.5cm}){%
     \parbox{0.62\textwidth}{\raggedright\fontsize{9}{11}\selectfont #1}}%
 }}%
\newcommand{\conf}[1]{%
\AddToShipoutPictureBG*{%
\AtPageUpperMyright{#1}
}
}

\makeatletter
\newcommand*{\rom}[1]{\expandafter\@slowromancap\romannumeral #1@}
\makeatother

\begin{document}

\maketitle
\conf{Proceedings of the 26th International Conference on Automation \& Computing, 
University of Portsmouth, Portsmouth, UK, 2-4 September 2021}
\thispagestyle{empty}
\pagestyle{empty}
\SetWatermarkText{Preprint}


\begin{abstract}
The main control tasks in autonomous vehicles are steering (lateral) and speed (longitudinal) control. PID controllers are widely used in the industry because of their simplicity and good performance, but they are difficult to tune and need additional adaptation to control nonlinear systems with varying parameters. In this paper, the longitudinal control task is addressed by implementing adaptive PID control using two different approaches: 
Genetic Algorithms (GA-PID) and then Neural Networks (NN-PID) respectively. The vehicle nonlinear longitudinal dynamics are modeled 
using Powertrain blockset library. Finally, simulations are performed to assess and compare the performance of the two controllers subject to external disturbances.

\end{abstract}

\section{INTRODUCTION}

The autonomous vehicle technology has been the main trend in the last decade. Although artificial intelligence and control researchers are quickly closing the gap to a real self-driving vehicle, several significant challenges still persist. A fully autonomous vehicle requires fully automatic steering and speed control modules. Cruise control systems (CC) control the speed of the vehicle and have been the subject of extensive research: Osman \textit{et al.} \cite{1} and Simorgh \textit{et al.} \cite{2} for instance showed that CC versions for highway driving worked only above certain speed limits as the braking part was not addressed \cite{3}. Adaptive CC systems were introduced later on to tackle both acceleration and braking \cite{4}, these systems are divided into high level and low level controllers. The first computes the required acceleration or deceleration, and the second actuates the throttle or the brake. These systems allowed to automate vehicles even for urban driving.

Many speed controllers are available in the literature. \cite{5} worked on the sliding mode control for vehicle platooning in uncertain topologies. \cite{6} designed a preview controller that integrates road slope, speed profile and vehicle dynamics to form an augmented optimal control. \cite{7} used the model predictive control (MPC) for low velocity tracking. Actually, most control systems are based on the PID controller which is easy to design and yet performs very well. However, the non-linearities of the vehicle subsystems and the variety of disturbances that face the vehicle such as road elevation, wind speed and road conditions require the PID to be adaptive to these changing parameters. \cite{2} and \cite{3} worked on adaptive PID control using the inverse model theory, \cite{8} used the twiddle algorithm to tune a PID for trajectory tracking and \cite{9} used Genetic Algorithms (GA) to tune a PID for controlling DC motors. \cite{10} and \cite{11} used NN-PID for frequency regulation and vehicle steering control respectively. In \cite{12}, the authors used radial basis function to adapt PID gains for the vehicle longitudinal control. The research of \cite{13} presented an adaptive PID with GA optimization for the lateral control task. However, most of these methods are based on oversimplified models excluding or simplifying external disturbances which significantly affects the performance in real-life situations. In addition, most adaptation methods require the measurement of certain parameters that are not easily accessible. To address this issue, two PID adaptation methods were applied to control the vehicle velocity; the first method is based on a GA with a complete vehicle model taking into account external disturbances, and the second one uses a feed forward neural network (NN) that utilises only the error and control signals to adapt the PID gains. Simulations were carried to assess the performance of the designed controllers. Section \rom{2} of this paper deals with the modeling of the vehicle longitudinal dynamics. The design, optimization and adaptation of the controller are presented in section \rom{3}. The performance of the controller is evaluated and discussed in section \rom{4}, and section \rom{5} concludes the paper and gives some perspectives for future work.
\vspace{-2pt}
\begin{figure}[!t]
\centering
\includegraphics[width=0.75\linewidth]{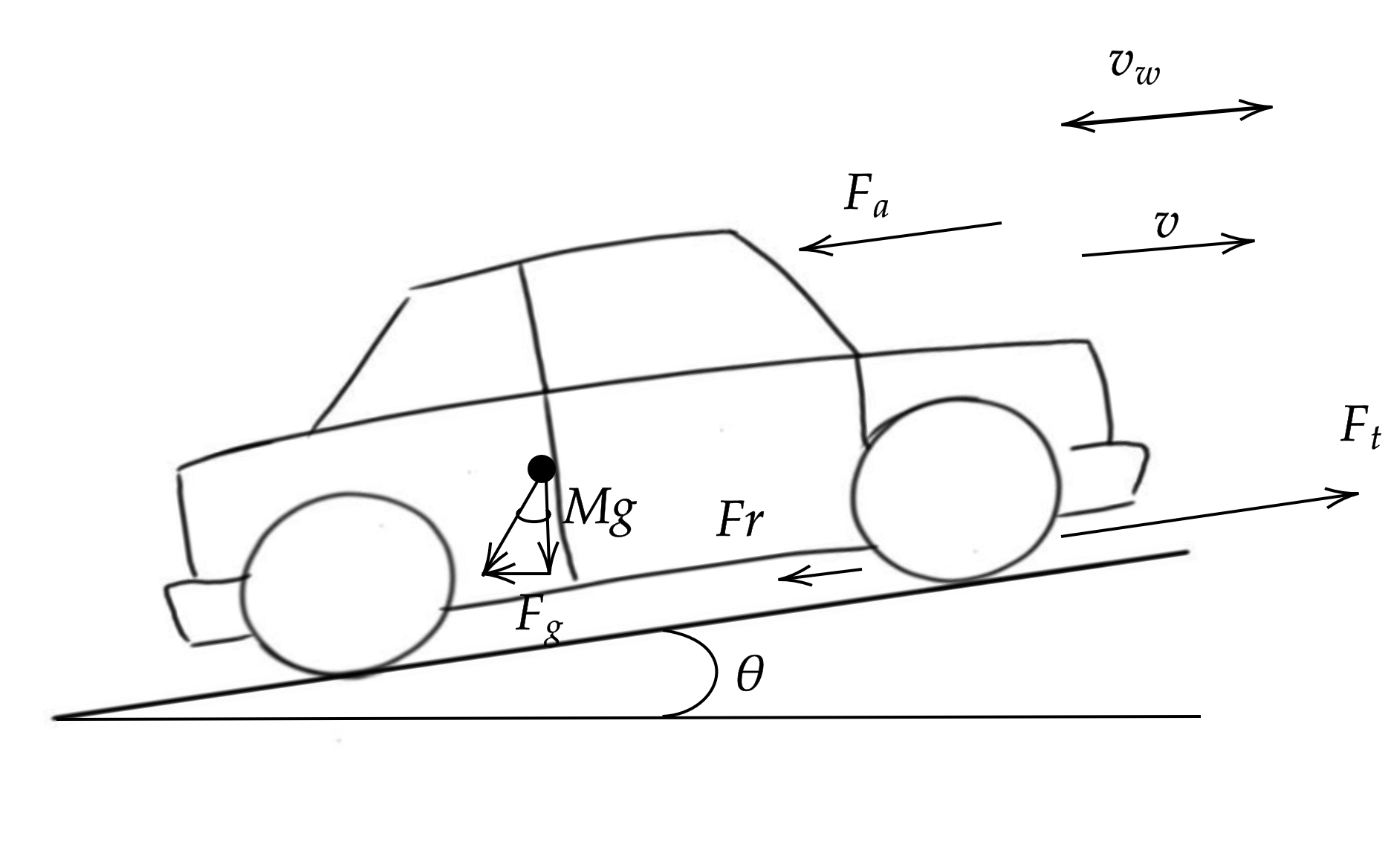}
\caption{Forces and disturbances affecting vehicle dynamics.}
\label{fig:1}
\end{figure}
\setlength{\textfloatsep}{5pt}
\section{Vehicle Modeling}\label{sec:Vehicle Modeling}

The longitudinal control of an autonomous vehicle regulates its speed to follow a speed profile. The latter is obtained from a speed planner or profiler based on changes of the external environment like the shape of the road, traffic signs and weather conditions. The vehicle is mainly propelled by the engine drive torque and stopped by the braking torque.

\subsection{Vehicle Body Dynamics}

The longitudinal dynamics of the vehicle are obtained by applying Newton’s second law as shown in Fig. \ref{fig:1}, this results in a simplified point mass model \cite{14}:

\begin{equation}
\label{eq29}
\left\{
    \begin{array}{ll}
        F_t &= M \frac{dv}{dt} + F_a + F_g +F_r \\
        F_g &= M g \sin \theta \\
        F_a &= \frac{1}{2} \rho C_d A v^2 \\
        F_r &= M g C_r \cos \theta
    \end{array}
\right.
\end{equation}
 
where the parameters are defined as:
\begin{itemize}
\item $M$   the total mass of the vehicle
\item $F_t$   the traction force (drive force)
\item $F_a$   the aerodynamic drag (air friction)
\item $F_g$   the downgrade force (gravity pull)
\item $F_r$   the rolling resistance force
\item $\rho$   the air density
\item $C_r$   the rolling resistance coefficient
\item $C_d$   the drag coefficient
\item $A$   the vehicle cross sectional area
\item $\frac{dv}{dt}$   the vehicle acceleration/deceleration.
\end{itemize}

\subsection{Engine Dynamics}

The engine and transmission dynamics make up the vehicle's propulsion system and are highly non-linear. For electric vehicles, the engine is an electric motor that supplies power to the wheels through the transmission. The latter is linked to the electric drive through a gearbox which can be modeled under the assumption of no losses by (\ref{eq2}): 

\begin{equation}
\label{eq2}
\left\{
    \begin{array}{ll}
        T_g &= k_g T_e \\
        \omega_g &= k_g \omega_e
    \end{array}
\right.
\end{equation}
where $\omega_g$ and $\omega_e$ are the angular speeds of the gear and the engine, $T_g$ and $T_e$ are the gear and electric drive torques and $k_g$ represents the ratio of the gearbox. The electric torque is generated by the electric machine and its model can be static or dynamic \cite{15}. The dynamic model depends on the type of machine whether it is a DC or AC (synchronous$\backslash$asynchronous) motor. It links the input current$\backslash$voltage to the output electric torque. The static model which is used in this work and given in (\ref{eq3}), is a data-driven model that is based on the electric drive efficiency map:

 \begin{equation}
\label{eq3}
\left\{
    \begin{array}{ll}
        T_e&= T_\text{e-ref}\\
        I_b&=\dfrac{T_e \omega_g}{u_b \eta^k}
    \end{array}
\right.
\end{equation}

 with $k=1$ if $T_e \omega_g \geq 0$, $k=-1$ otherwise. $T_\text{e-ref}$ represents the reference electric torque obtained from the motor efficiency map, $u_b$ denotes the battery voltage, $I_b$ and $\eta$ are the battery current and motor efficiency respectively. Positive currents discharge the battery while negative currents, as in regenerative braking, charge the battery.

\subsection{Wheel Dynamics}
The wheel translates the gear torque into a traction force that propels the vehicle. Thus, the angular speed of the wheel determines the velocity of the vehicle as the following:
\begin{equation}
\label{eq4}
\left\{
    \begin{array}{ll}
        F_{t/b}&= \frac{1}{R_w} T_{g/b}\\
        \omega_w&= \frac{1}{R_w} V
    \end{array}
\right.
\end{equation}

\noindent where $F_{t/b}$ and $T_{g/b}$ denote either the tractive force (Fig. \ref{fig:1}) and the drive torque or the braking force and braking torque respectively. $R_w$ and $V$ are the respective wheel radius and vehicle longitudinal velocity.  

\subsection{Tire Dynamics}
The tires of a vehicle undergo different kinds of deformations and behave in a highly nonlinear manner, the most widely used longitudinal tire model is the magic formula or Pacejka model \cite{16}:
\begin{equation}
\label{eq5}
F_x = F_z D \sin(C \tan^{-1}[{Bk-E[Bk-\tan^{-1}(Bk)]}])
\end{equation}
\noindent with $F_x$ and $F_z$ being the tire longitudinal force and vertical load, $k$ representing the wheel slip and $B,C,D$ and $E$ being empirical tire coefficients.

\subsection{Brake Dynamics}
The brakes of the vehicle can be either a drum or a disk system. The disk braking system has been used in this paper for both front and rear wheels. The braking torque delivered by the mechanism is given as:
\begin{equation}
\label{eq6}
T_b= \frac{\mu P \pi B_a^2 R_m N_{pads}}{4}
\end{equation}
\noindent where $P, R_m$ and $B_a$ are the applied brake pressure, the mean radius of the brake pad and the brake actuator diameter respectively. $ \mu$ and $N_{pads}$ are the disc-pad friction coefficient and the number of pads respectively.

\begin{figure*}[!h]
\centering
\includegraphics[width=0.8\linewidth]{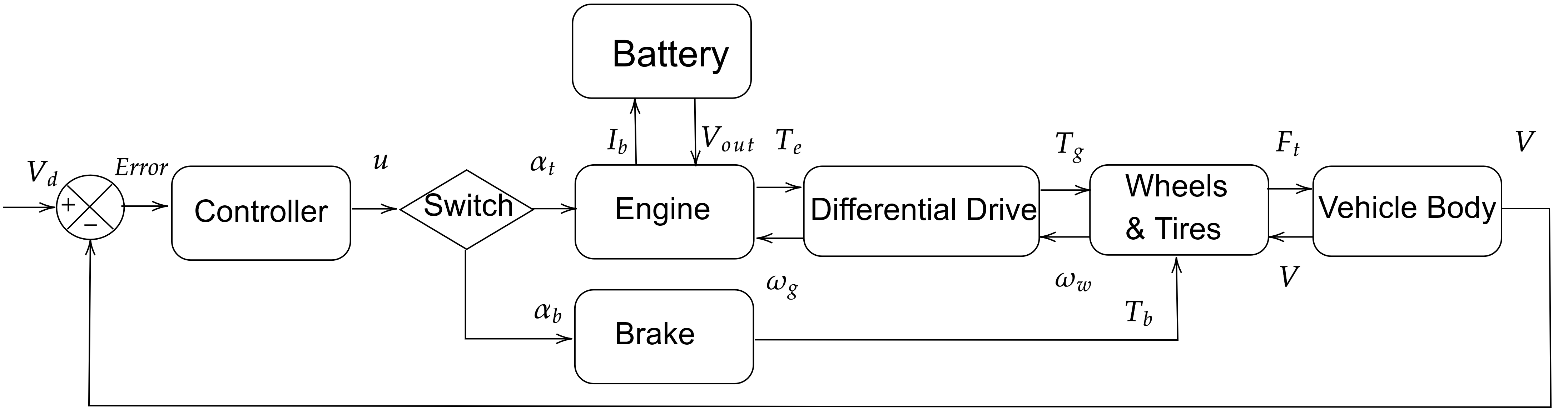}
\caption{Model Block Diagram}
\label{fig:2}
\end{figure*}

\subsection{Battery modeling}
The electric motor is powered by the battery voltage $u_b$ which can be modeled by a resistor $R_b$ in series with an open circuit voltage $V_{oc}$. Both of them depend on the battery state of charge (SOC):
\begin{equation}
\label{eq7}
\left\{
\begin{array}{ll}
  u_b&=  N_s(V_{oc}(SOC) + \frac{I_b}{N_p} R_b(SOC))\\
  SOC &= \frac{1}{C_b} \int^t_0 \frac{I_b}{N_p} dt
\end{array}
\right.
\end{equation}
\noindent where $N_s$, $N_p$ are the number of cells in series and in parallel respectively, and the denominator $C_b$ is the battery capacity.

The Powertrain blockset of MATLAB has been used to build a high fidelity vehicle longitudinal model (see Fig. \ref{fig:2}). The vehicle, motor and battery characteristic parameters were based on the Renault Zoe Vehicle (see table \ref{tab:4}). The rest of the parameters are blockest default values. Moreover, both the brake and accelerator pedal actuators are modeled as a first order system with time constants $\tau_a=0.75s $ and $\tau_b=1s $ respectively \cite{4}. 
The controller generates either a braking or an acceleration command at a time. Based on the sign of the control signal, the switching logic sends the command to the braking system or to the engine such that a negative control signal means that the braking is activated, otherwise, acceleration is activated \cite{12}.

\begin{table}[!h]
\caption{Simulation parameters.} 
\label{tab:4}
\centering
\begin{tabular}{c| c c| c} 

$M$ & $1468kg$ & $C_d$ & $0.29$\\ [0.5ex] 

$C_r$ & $0.007$ & $k_g$ & $3.4$\\

$A$ & $2.22m^2$ & $R_w$ & $0.329m$\\

$\rho$ & $1.225kg/m^3$ & $N_s$ & $96$ \\

$\mu$  & $0.9$ & $g$ & $9.81m/s^2$\\

CG height  & $0.35m$ & $N_{pads}$ & $2$ \\

Front axle to CG & $1.455m$ & $\eta$  & $95\%$\\

Rear axle to CG & $1.132m$& $C_b$ & $132Ah$\\

$R_m$ & $0.1778m$ & $N_b$ & $2$\\[1ex] 

\end{tabular}
\end{table}

\vspace{-10pt}
\section{Controller Design}

\subsection{PID Controller}

The PID controller is heavily used in the industry due to its simple design and satisfactory performance. The input to the controller is the error signal defined as the difference between the set-point and the system output as in Fig. \ref{fig:7}. The output of the controller depends on the PID gains and is governed by the following equation:
\begin{equation}
\label{eq8}
\begin{split}
u(t) = K_p e(t) + K_i \int e(\tau)d\tau + K_d \frac{de(t)}{dt}
\end{split}
\end{equation}

\noindent where $K_p$, $K_i$ and $K_d$ are the proportional, integral and derivative gains and $e(t)$ is the error between set-point $R$ and output $Y$. The discrete form of the PID has been used in this work where the control law is obtained using the trapezoid method as follows:

\begin{equation}
\label{eq9}
\left\{
\begin{array}{ll}
\Delta u(k) & =  \alpha e(t) + \beta e(k-1) + \gamma e(k-2)\\
\alpha & = K_p + \frac{K_d}{T_s} + \frac{K_i T_s}{2}\\
\beta & = \frac{K_i T_s}{2} -2 \frac{K_d}{T_s} -K_p\\
\gamma & = \frac{K_d}{T_s}
\end{array}
\right.
\end{equation}
with $T_s$ being the sampling time.\\

The performance of a PID controller is measured in terms of overshoot $O_s$, rising time $R_t$ and settling time $S_t$. It depends on the choice of gains which is not straight forward, there are several tuning techniques such as root locus, gain-phase margin and the widely used Ziegler Nichols \cite{17}. However, these methods are generally limited to linear systems and have poor disturbance rejection. To account for the effects of road angle and wind speed on the controller performance, the designed model is simulated with GA to optimize the gains of the PID controller. The choice of this meta-heuristic optimization method is based on its good performance reported in the literature \cite{9}.

\subsection{GA-based PID Optimization}

\begin{figure}[b]
\centering
\includegraphics[width=0.45\textwidth]{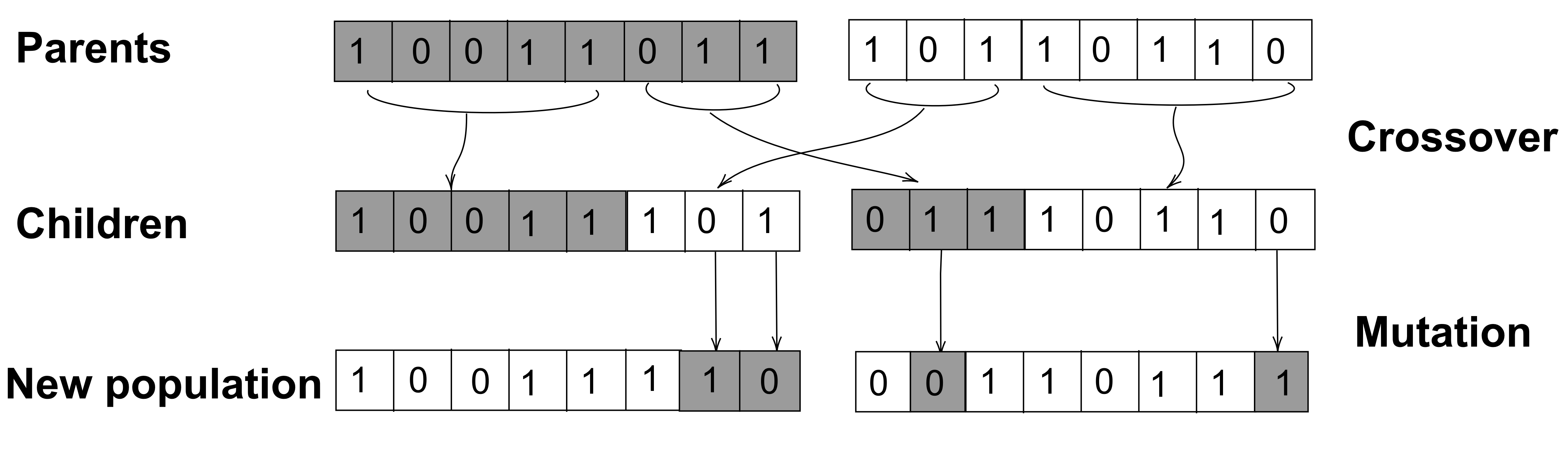}
\caption{Genetic operations}
\label{fig:4}
\end{figure}

Genetic algorithms are inspired from biological evolution theory which is based on the process of natural selection and genetics \cite{9}. The algorithm uses population sets of multiple potential encoded solutions called chromosomes. Each chromosome is rated by a fitness function (\ref{eq10}) that determines how good of a solution it is.
\begin{equation}
\label{eq10}
\text{Fitness}_\text{value} = \frac{1}{\text{Performance}_\text{index}}
\end{equation}
GA algorithms use crossover and mutation operations (Fig. \ref{fig:4}) to enhance the existing solutions. 

The crossover operation merges existing chromosomes together, while the mutation operation modifies the encoding of existing chromosomes to yield new and better chromosomes. The fitness score of each chromosome is assigned by an objective function and only fit chromosomes are selected to evolve in the next operation. In the context of PID tuning, the objective function minimizes the error that is fed to the controller. The mean of squared error (MSE) has been used as the objective function for a consistent performance comparison between GA-PID and NN-PID. 

\begin{figure}[t]
\centering
\includegraphics[width=0.35\textwidth]{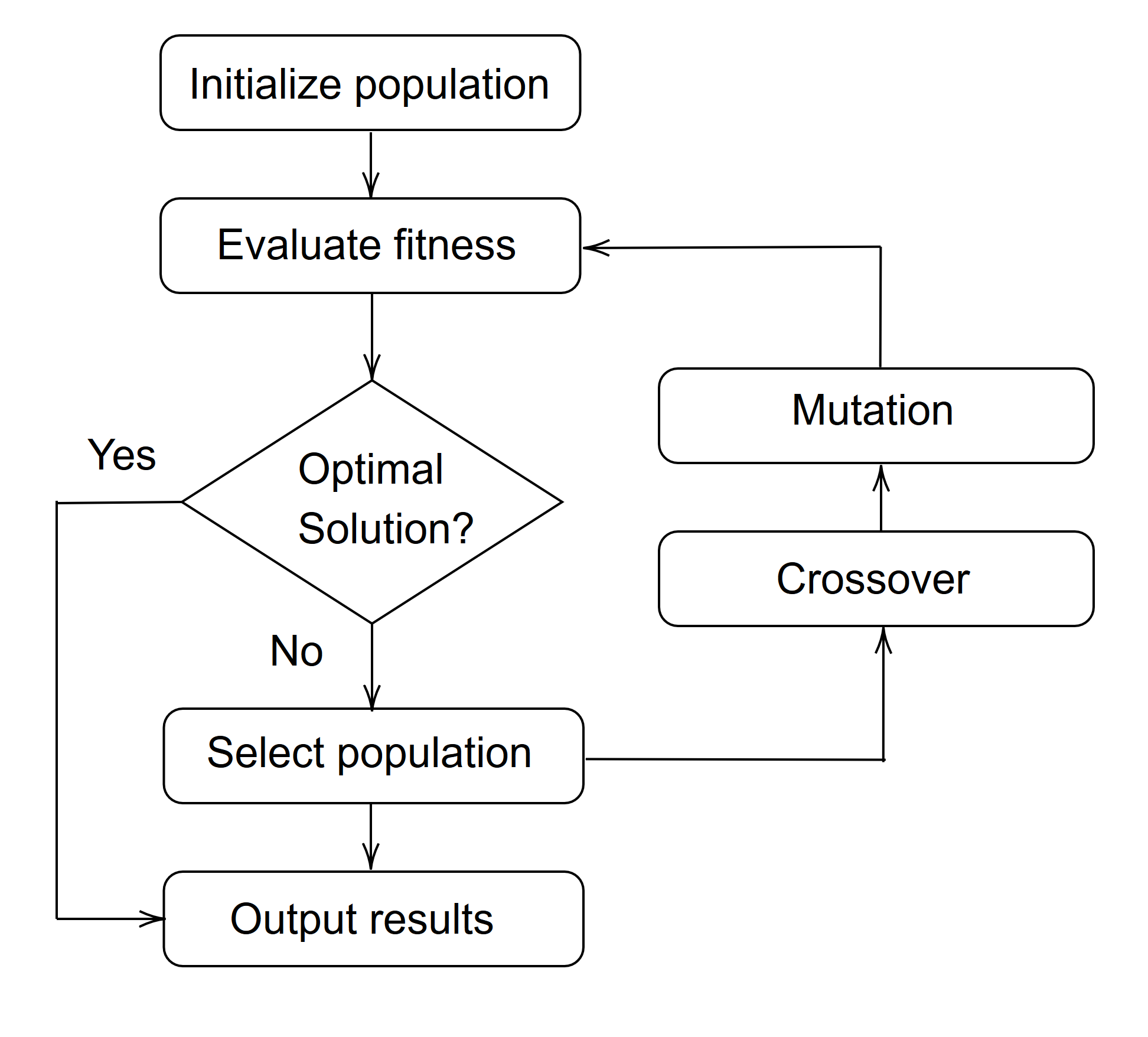}
\caption{Genetic algorithm process.}
\label{fig:5}
\end{figure}

The GA optimization tool in (Fig.  \ref{fig:5}) has been used. The process initiates a random population of a certain size containing values (chromosomes) for the PID gains ($K_p$, $K_i$ and $K_d$) with which the vehicle longitudinal model is simulated. 
\setlength{\textfloatsep}{5pt}
The fitness value of each chromosome is assessed through the MSE. Fit chromosomes with high fitness scores are selected as fit parents, while unfit chromosomes with low fitness scores are discarded. The GA then performs mutation and crossover on the best fit parents and produces a new population. The number of best fit parents to be enhanced is predefined as a percentage called elite count value. This whole process is repeated over many generations. Moreover, the first selection of population individuals is random but bounded between a lower and an upper limit.
The PID gains are optimized with regard to the varying parameters of the model and new optimal PID gains are generated for each case of the varying parameters. A variety of cases were tested and the results were stored in a small data-set. The data-set is then used in a lookup table to select the optimized gains that correspond to the varying parameters. This method allows continuous adaptation of the PID gains as the varying parameters change along a speed profile or drive cycle. 

\subsection{Neural Network PID Adaptation}
\begin{figure}[b]
\centering
\includegraphics[width=0.35\textwidth]{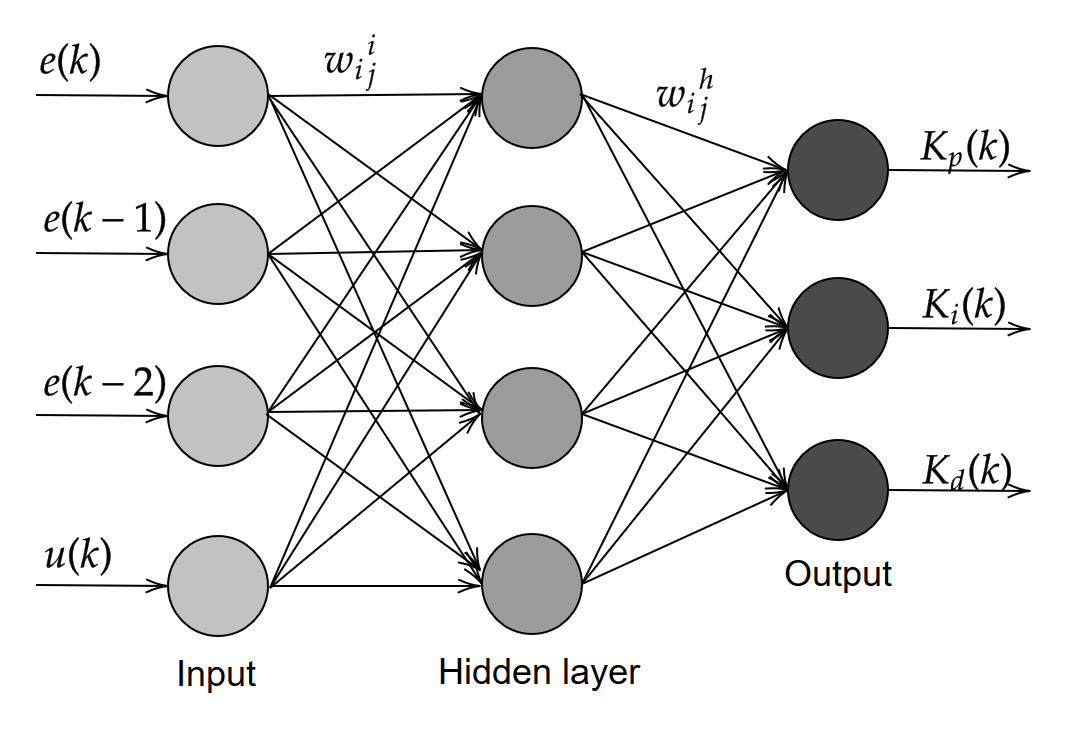}
\caption{Neural network architecture.}
\label{fig:6}
\end{figure}

\begin{figure}[t]
\centering
\includegraphics[width=0.4\textwidth]{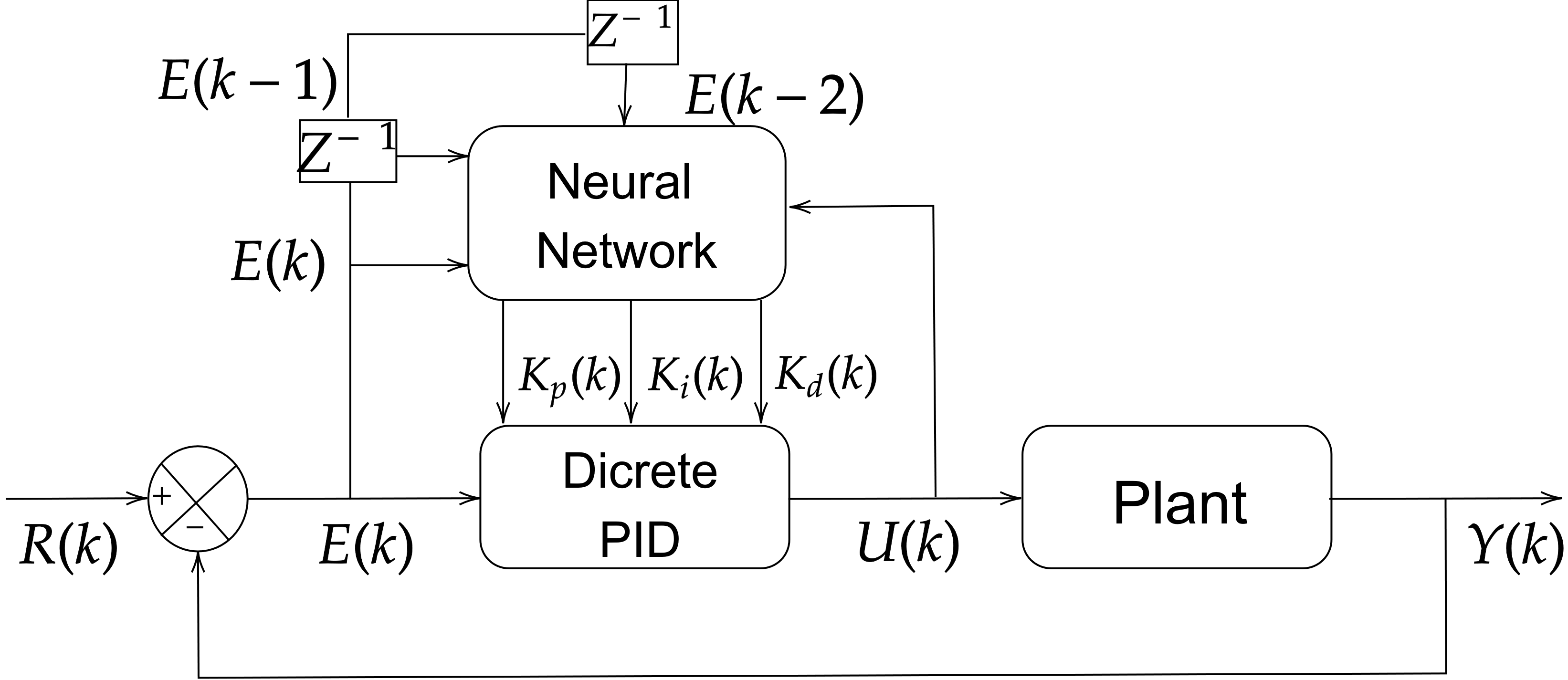}
\caption{Block diagram of the system}
\label{fig:7}
\end{figure}

The use of a NN inside the control loop makes it possible to adapt the PID gains online. This method presents a better and more suitable alternative to offline optimization with GA. 
The NN is trained online using the back-propagation of the error generated from previous PID gains. The architecture of the NN consists of an input layer with four inputs, a hidden layer with four neurons and an output layer with three neurons (Fig. \ref{fig:6}).
The sigmoid activation function has been used for the hidden layer, while contrary to existing works, instead of the identity, the ReLU function has been used for the output layer to avoid predicting negative gains. The overall system is given in Fig. \ref{fig:7}.
\setlength{\textfloatsep}{5pt}

The inputs to the hidden layer $(S_j)$ and output layer $(Z_j)$ are given by (\ref{eq13}). $w_{ij}^{i}$ are the weights between input and hidden layer and $w_{ij}^h$ are the weights between hidden and output layer, $n_h$ and $n_o$ represent the number of hidden units and output units respectively. $x_i$ are the inputs to the network, the outputs of the hidden layer $h_j$ and the outputs of the network $o_j$ are given in (\ref{eq14}) where $f$ and $g$ are the sigmoid and ReLU
activation functions respectively.   
\begin{equation}
\label{eq13}
\left\{
    \begin{array}{ll}
        S_j =& \sum_{i=1}^{n+1}w_{ij}^i x_i,  \; \text{for} \; j=1,..,n_h\\
        Z_j =& \sum_{i=1}^{n_h+1}w_{ij}^h h_i,  \; \text{for} \; j=1,..,n_0
    \end{array}
\right.
\end{equation}
\begin{equation}
\label{eq14}
\left\{
    \begin{array}{ll}
        h_j =& f(S_j)\\
        o_j =& g(Z_j)
    \end{array}
\right.
\end{equation}
\setlength{\textfloatsep}{5pt}
\setlength{\textfloatsep}{5pt}
The gradient descent algorithm minimizes the MSE function  $E(k)=\frac{1}{2}e(k)^2$, which is back-propagated through the network using the chain rule \cite{13}. The gradient of the error term is given by:
\begin{equation}
\label{eq15}
\frac{\delta E(k)}{\delta w_{ij}^h(k)}=\frac{\delta E(k)}{\delta e(k)} \frac{\delta e(k)}{\delta y(k)} \frac{\delta y(k)}{\delta u(k)} \frac{\delta u(k)}{\delta o_j} \frac{\delta o_j}{\delta w_{ij}^h}
\end{equation}
where
\begin{equation}
\label{eq16}
\left\{
    \begin{array}{ll}
        \frac{\delta E(k)}{\delta e(k)}& \frac{\delta e(k)}{\delta y(k)} = e(k)\\
        \frac{\delta y(k)}{\delta u(k)} &\approx sign(\frac{Y(k)-Y(k-1)}{U(k)-U(k-1)})\\
        \frac{\delta u(k)}{\delta o_1} &= e(k) - e(k-1)\\
        \frac{\delta u(k)}{\delta o_2} &= \frac{T_s}{2}(e(k)+e(k-1))\\
        \frac{\delta u(k)}{\delta o_3} &= \frac{1}{T_s}(e(k)-2e(k-1)+e(k-2))\\
        \frac{\delta o_j}{\delta w_{ij}^h} &= h_i\\
    \end{array}
\right.
\end{equation}

The term $\frac{\delta y(k)}{\delta u(k)}$ is not accessible in the network, therefore, it is approximated by sign$(\frac{Y(k)-Y(k-1)}{U(k)-U(k-1)})$, which shows whether the control law is improving. To avoid singularities, the expression sign{\small $((Y(k)-Y(k-1))\times(U(k)-U(k-1))$} is used instead. For the hidden layer, the error term is defined as:
\begin{equation}
\label{eq17}
\frac{\delta E(k)}{\delta w_{ij}^i(k)} = \frac{\delta E(k)}{\delta h_j(k)} \frac{\delta h_j(k)}{\delta S_j(k)} \frac{\delta S_j(k)}{\delta w_{ij}^i(k)}
\end{equation}
where: $i = 1,..,n_i+1$, and $j=1,..,n_h$. 
\begin{equation}
\label{eq18}
\left\{
    \begin{array}{ll}
         \frac{\delta E(k)}{\delta h_j(k)} &= -e(k) \frac{\delta y(k)}{\delta u(k)} \sum_{i=1}^{n_0}  \frac{\delta u(k)}{\delta o_l(k)} w_{jl}^h(k)\\
         \frac{\delta h_j(k)}{\delta S_j(k)} &= h_j(k)(1-h_j(k))\\
         \frac{\delta S_j(k)}{\delta w_{ij}^i(k)} &= x_i(k)
        
    \end{array}
\right.
\end{equation}
The weights are updated according to the gradient descent rule as follows:
\begin{equation}
\label{eq19}
\left\{
    \begin{array}{ll}
         w_{ij}^h(k+1) &= w_{ij}^h(k)- \eta \frac{\delta E(k)}{\delta w_{ij}^h(k)}\\
              w_{ij}^i(k+1) &= w_{ij}^i(k)- \eta \frac{\delta E(k)}{\delta w_{ij}^i(k)}
    \end{array}
\right.
\end{equation}
with $\eta$ being the learning rate.

\section{Simulation Results and Discussion}

\begin{figure}[hbt!]
\centering
\includegraphics[width=0.45\textwidth]{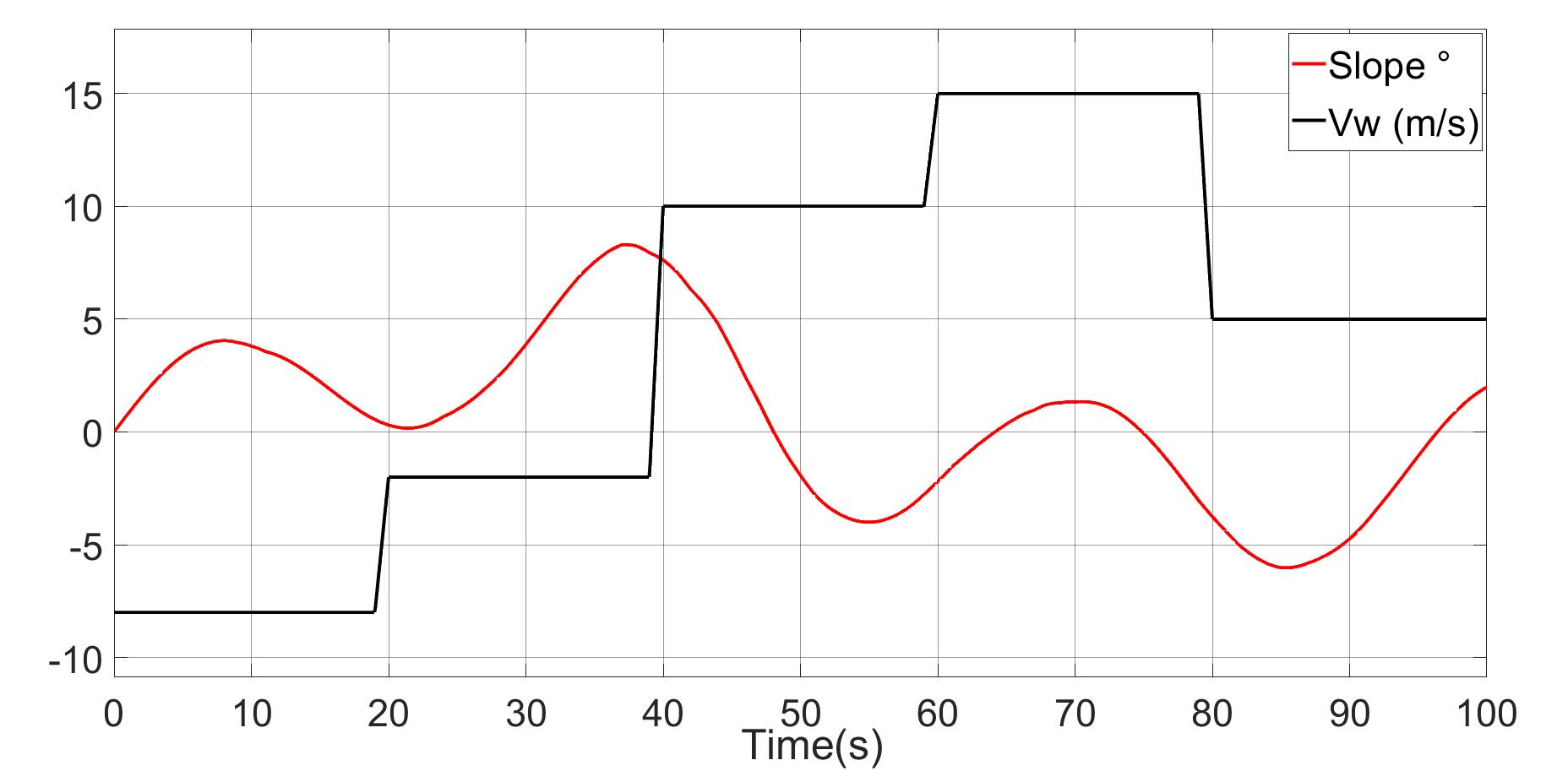}
\caption{Wind speed and road slope profiles}
\label{fig:8}
\end{figure}

The designed PID controllers have been tested  with the vehicle longitudinal model given in section \ref{sec:Vehicle Modeling}. The GA algorithm is used to optimize the gains of the PID (GA-PID) for three scenarios:
\begin{itemize}
\item \textbf{Sc1} A varying speed profile without disturbances ($V_w=0$ and $\theta=0$).
\item \textbf{Sc2} A varying speed profile subject to varying wind speed and road slope (see Fig. \ref{fig:8}).
\item \textbf{Sc3} Optimization for different values of $V_\text{ref}$,$V_w$ and $\theta$ separately.
\end{itemize}
\setlength{\textfloatsep}{2pt}
\setlength{\textfloatsep}{2pt}
Both scenarios 1 and 2 produce a triplet of optimal gains ($K_p, K_i$ and $K_d$) for the specific speed and disturbance profiles. Scenario 2 shows the effect of the disturbances ($v_w$ and $\theta$) on the controller performance. Scenario 3 is adaptive and more generalized and produces a data-set of optimized gains corresponding to the changing values of $V_\text{ref},V_w,$ and $\theta$. 
The GA-PID is compared to the NN-PID with optimized initial weights and learning rate. Initialization of the weights and choice of the learning rate are very important stages for the performance of the NN-PID and its convergence speed. The learning rate is selected by trial and error after the first few simulations. 

\begin{figure}[hbt!]
\centering
\includegraphics[width=0.44\textwidth]{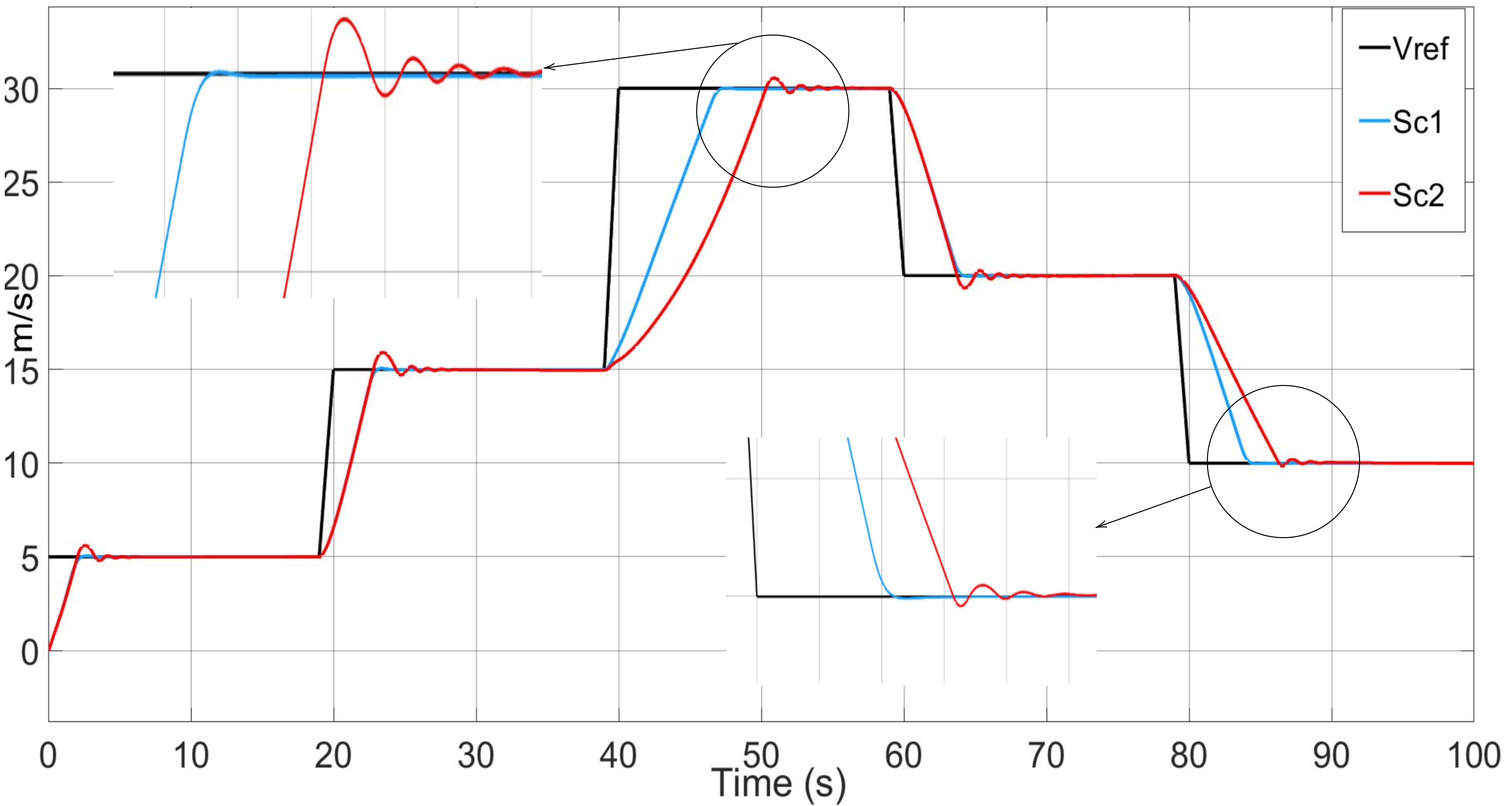}
\caption{Speed tracking for $Sc_1$ and $Sc_2$}
\label{fig:9}
\end{figure}

\begin{figure}[hbt!]
\centering
\includegraphics[width=0.42\textwidth]{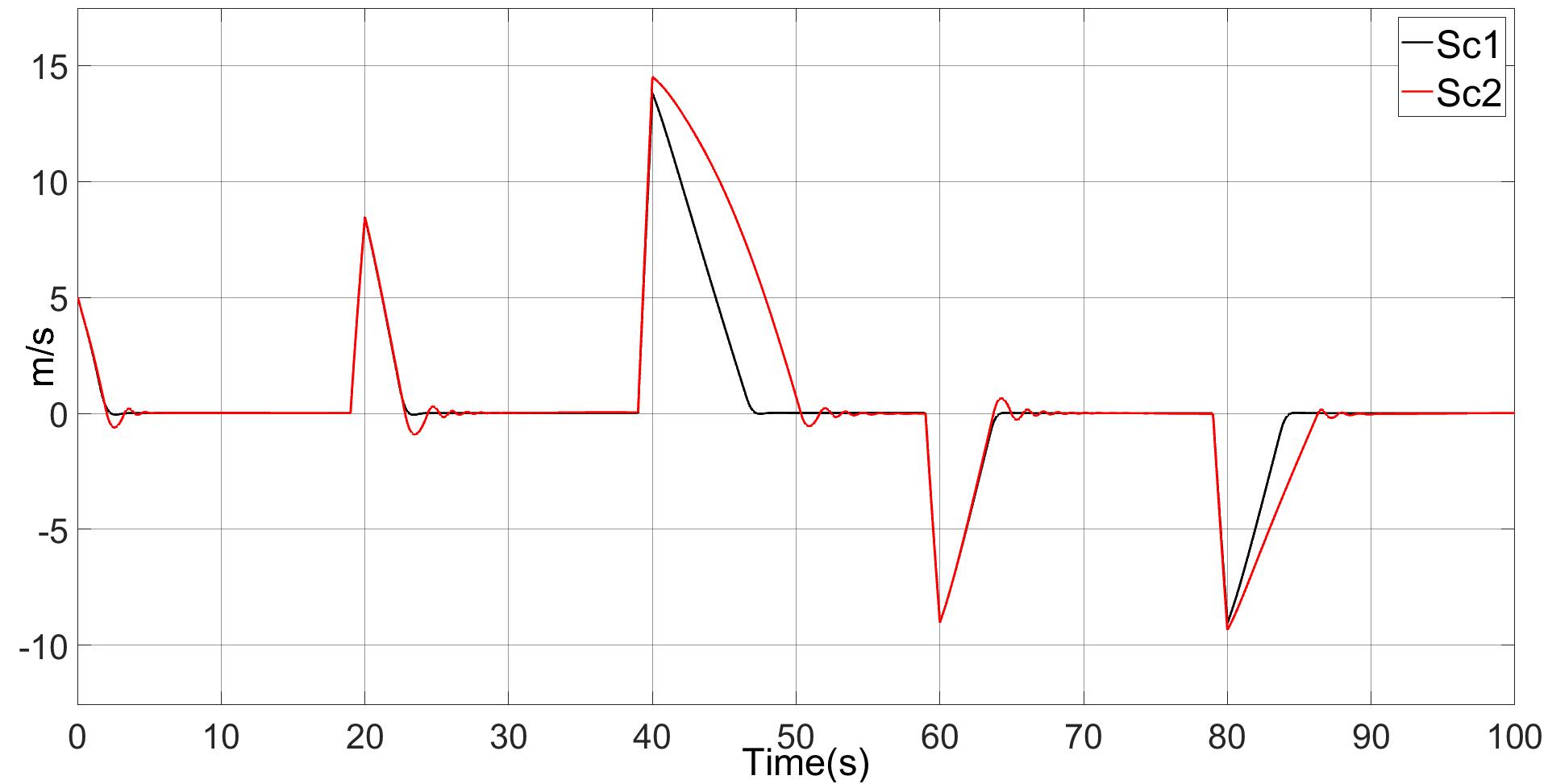}
\caption{Error variation for $Sc_1$ and $Sc_2$}
\label{fig:10}
\end{figure}

\begin{figure}[hbt!]
\centering
\includegraphics[width=0.42\textwidth]{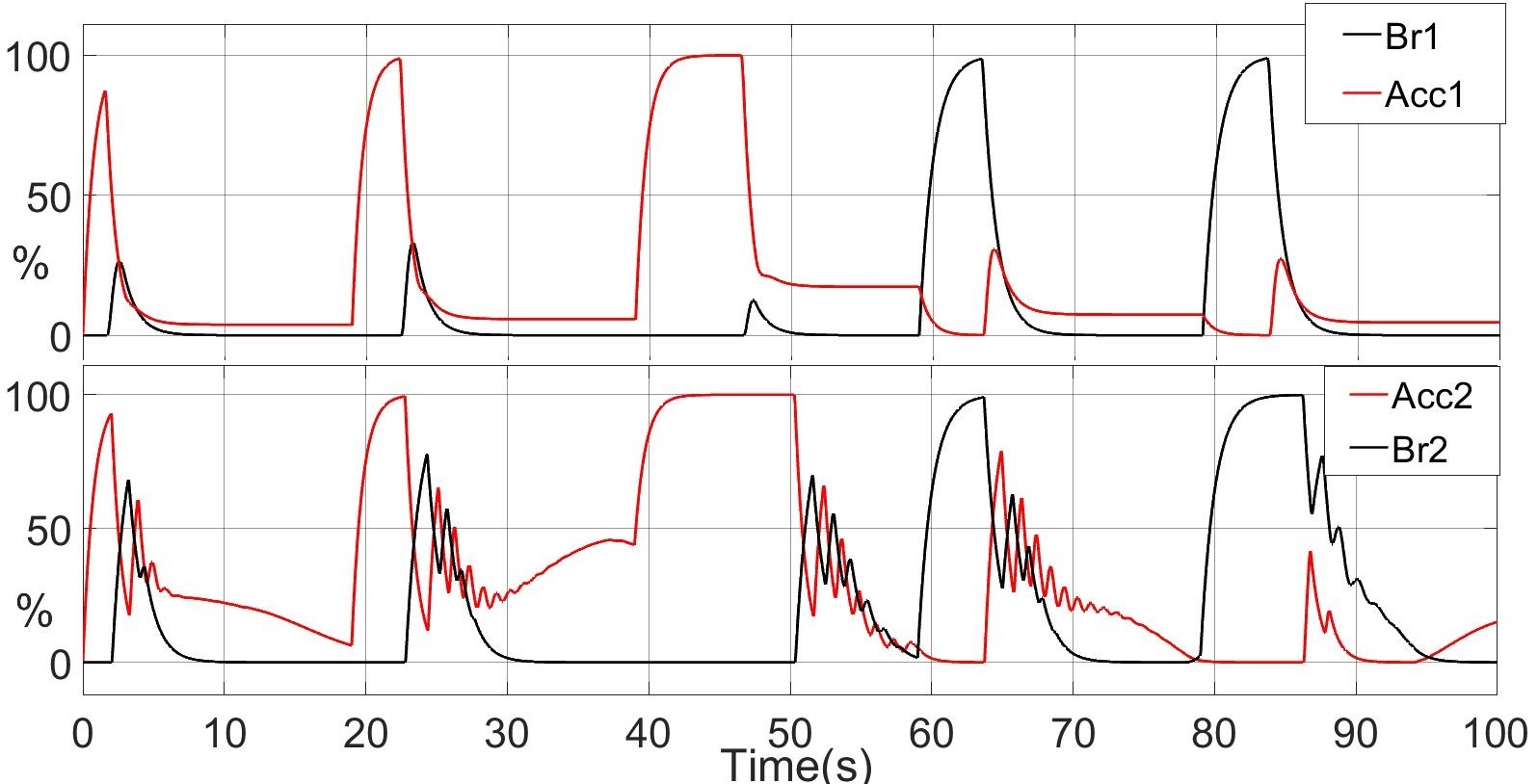}
\caption{Acceleration/Brake variation for $Sc_1$ and $Sc_2$}
\label{fig:11}
\end{figure}

\begin{figure}[hbt!]
\centering
\includegraphics[width=0.46\textwidth]{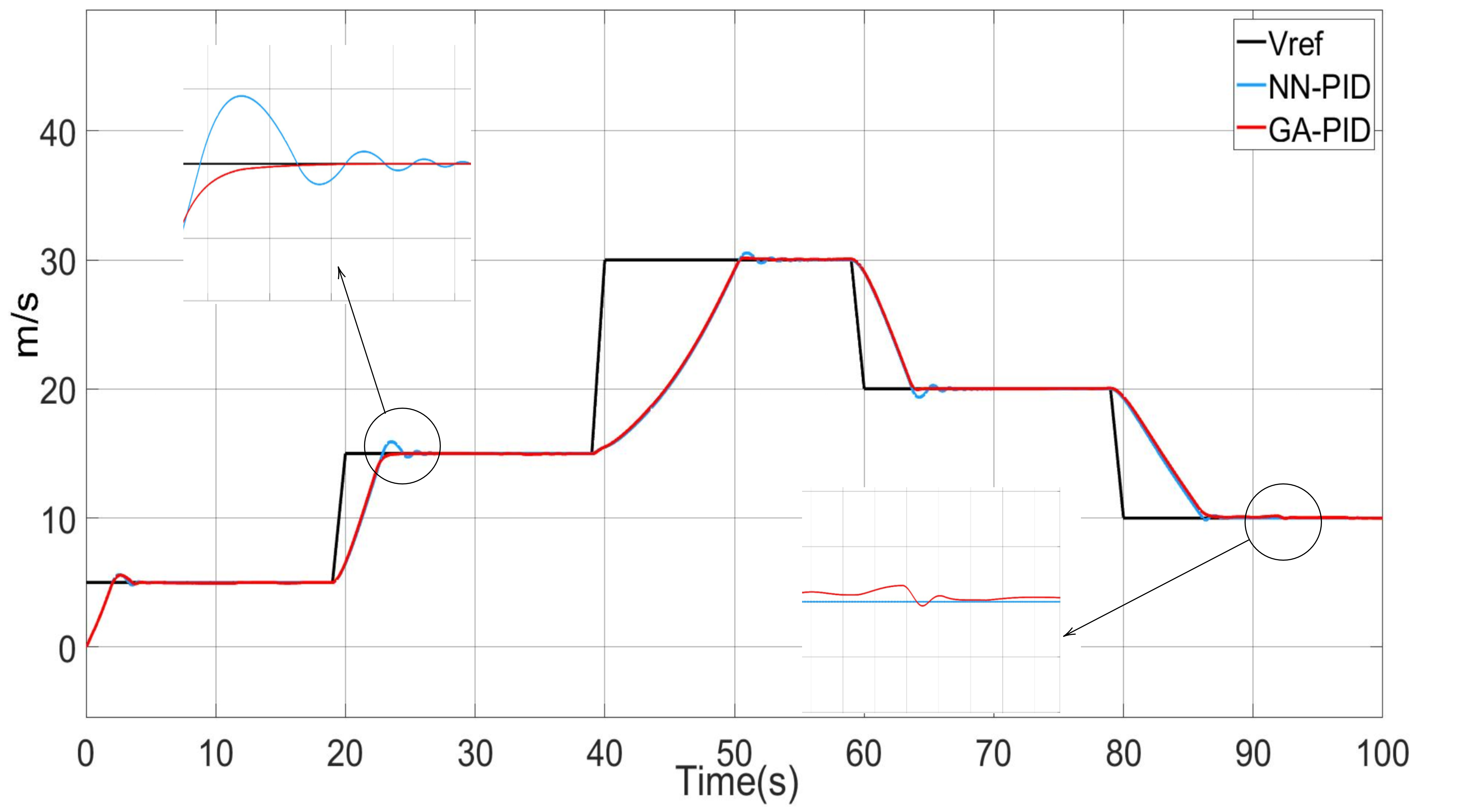}
\caption{Speed tracking for $Sc_3$ vs NN-PID.}
\label{fig:12}
\end{figure}

\begin{figure}[hbt!]
\centering
\includegraphics[width=0.42\textwidth]{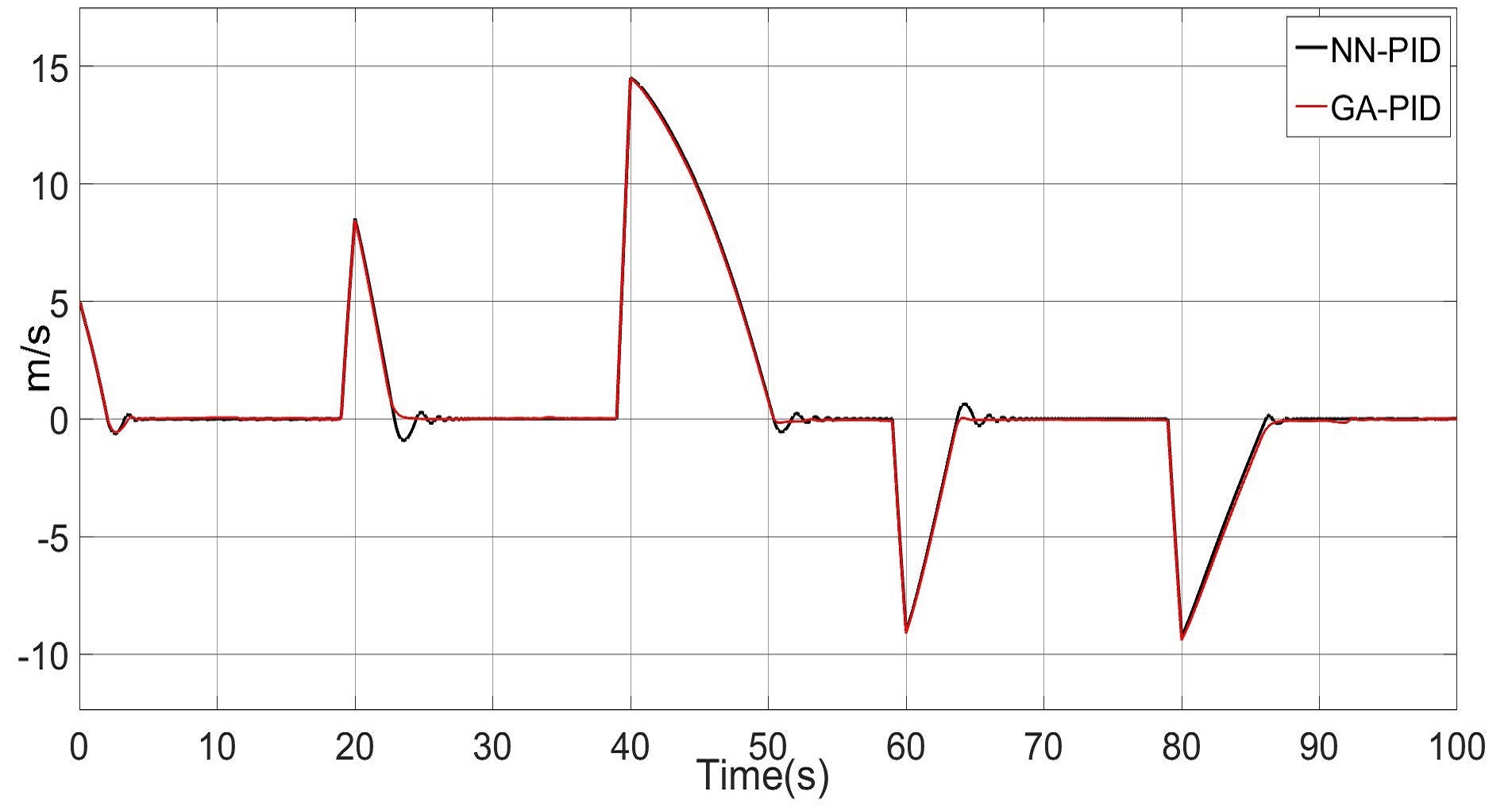}
\caption{Error variation for $Sc_3$ vs NN-PID.}
\label{fig:13}
\end{figure}

\begin{figure}[hbt!]
\centering
\includegraphics[width=0.42\textwidth]{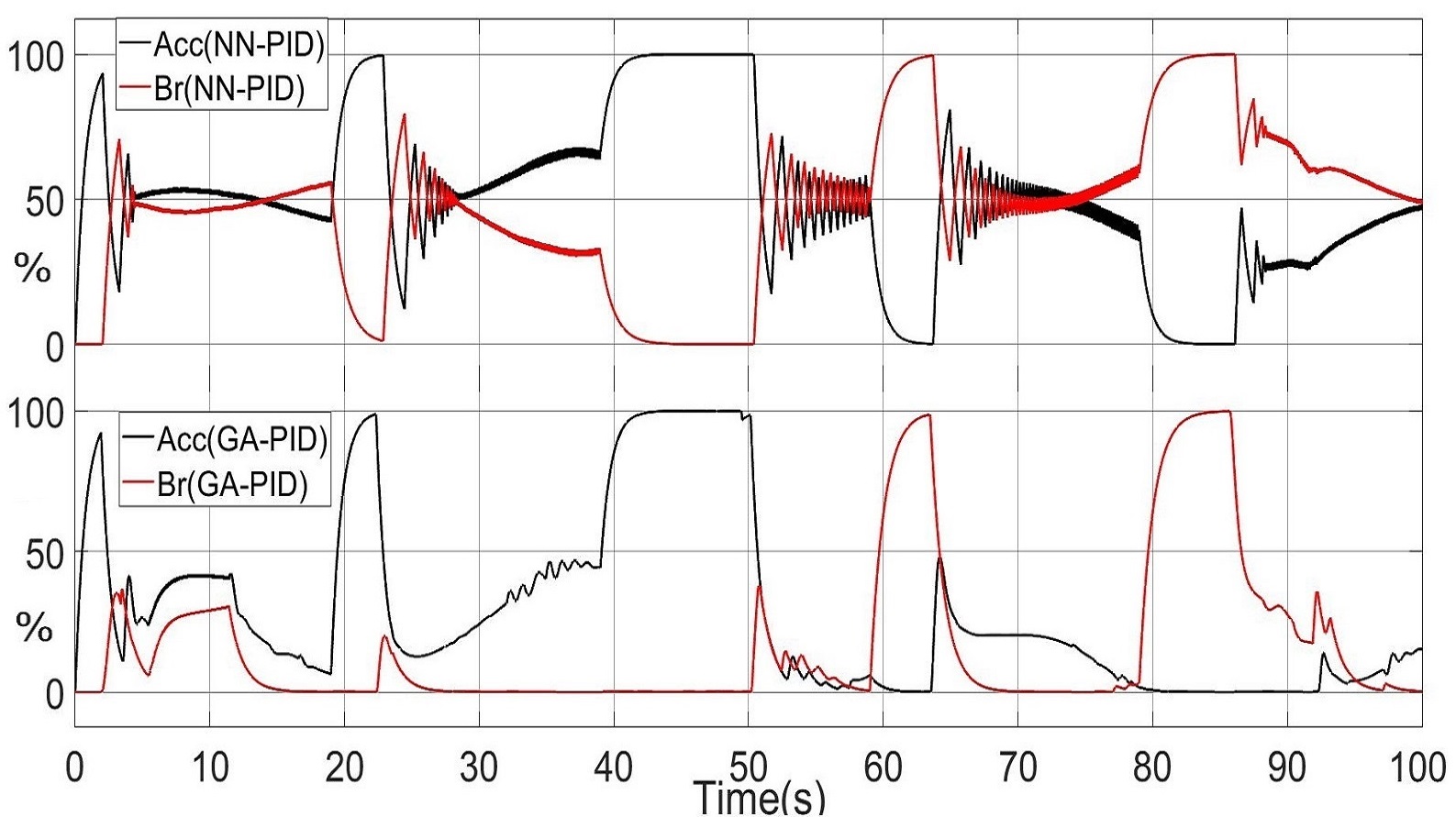}
\caption{Acceleration/Brake variation for $Sc_3$ vs NN-PID.}
\label{fig:14}
\end{figure}

The GA optimization for scenarios 1 and 2 results in the triplets $\{K_{p_1}=999.75$, $K_{i_1}=0.1$, $K_{d_1}=0.3\}$ and $\{K_{p_2}=383.79$, $K_{i_2}=0.11$, $K_{d_2}=128.84\}$ respectively. The corresponding results are compared in Fig. \ref{fig:9}-\ref{fig:11} where the speed tracking and error variation signals are in (m/s) and the acceleration and braking signals are in percentage (\%). 
For scenario 3, the GA-PID gains are optimized for $V_\text{ref}\in[0,30]$ m/s, $\theta \in[-10^{\circ},10^{\circ}]$ and $V_w \in[-10,15]$ m/s. The performance of PID adaptation using the generated data-set in a lookup table is compared to that of the NN-PID in Fig. \ref{fig:12}-\ref{fig:14}.

In scenario 1, the GA-PID produces a smoother and faster response with no overshoots due to the absence of disturbances. On the other hand, the inclusion of road slope and wind speed profiles in scenario 2 increases the response time which corresponds to real-life situations, this results in slower tracking with slight overshoots resulting in less smooth responses. 

\setlength{\textfloatsep}{2pt}
\setlength{\textfloatsep}{2pt}
\setlength{\textfloatsep}{2pt}
\setlength{\textfloatsep}{2pt}

The GA-PID in scenario 3 is smoother than the NN-PID which is more aggressive since the adaptation takes place at each iteration. On the other hand, the NN-PID is more robust to disturbances and faster because it reacts immediately whenever the error magnitude changes and adapts the PID gains to quickly reject the disturbances. Overall, both GA-PID and NN-PID insure adaptive control with very good performance and disturbance rejection.

\setlength{\textfloatsep}{2pt}
\setlength{\textfloatsep}{2pt}
\newpage
Table \ref{tab:3} compares the controller performance regarding the $MSE$ error and evaluates $R_t$, $S_t$ and $O_s$ for a step response using different values for the hidden layer neurons (h) and learning rate (lr). The GA-PID in scenario 2 has the best performance, while the NN-PID with (h=10, lr=0.01) is the second best controller.

\begin{table}[htb!]
\caption{Controller performance for different $h$ and lr.} 
\label{tab:3}
\centering
\begin{tabular}{c c c c c} 

\textbf{Method} & \textbf{MSE} & $\boldsymbol{R}_t(s)$& $\boldsymbol S_t(s)$& $\boldsymbol O_s(\%)$  \\ [0.5ex] 
\hline
NN-PID (lr=0.01, $h=4$) & 7.596 & 0.348 & 0.445 & 4.6\\[0.8ex] 

NN-PID (lr=0.1 , $h=4$) & 7.605 & \textbf{0.37} & 0.445 & 4.59\\[0.8ex] 

NN-PID (lr=0.8 , $h=4$) & 7.604 & \textbf{0.37} & 0.445 & 4.59\\[0.8ex] 

NN-PID (lr=0.01 , $h=10$) & 7.555 & 0.348 & 0.445 & 4.6\\[0.8ex] 

NN-PID (lr=0.8 , $h=15$) & 7.597 & 0.348 & 0.445 & 4.59\\[0.8ex] 

NN-PID (lr=0.8 , $h=20$) & 7.557 & 0.348 & 0.445 & 4.6\\[0.8ex] 

GA-PID(Sc1) & 7.563 & 0.87 & 1.164 & 6.68 \\[0.8ex] 

GA-PID(Sc2) & \textbf{7.502} & 0.348 & \textbf{0.444} & 4.79 \\[0.8ex] 

GA-PID(Sc3) & 7.534& 0.449& 0.616& \textbf{1.52} \\ [0.8ex] 

\end{tabular}
\label{table:nonlin} 
\end{table}
\vspace{-5pt}
\section{CONCLUSIONS}
This paper dealt with the longitudinal control of autonomous vehicles where the vehicle dynamics were modeled using the Powertrain blockset. The adaptive PID controller was designed using two different approaches. The first approach (GA-PID) uses GA where the MSE has been used as a cost function. The advantages of this approach are \begin{enumerate*}[label=\textit{(\alph*)}]
  \item  optimization with a complete vehicle model unlike other approaches which use simplified models in the form of transfer function or state space representation.
  \item optimization for a variety of cases and adaptation via a look-up table.
\end {enumerate*}
The downside of this approach is the heavy calculations required. The second approach (NN-PID) was based on online learning with an MLP NN where the PID gains are adapted with back-propagation. This approach has the following advantages:
\begin{enumerate*}[label=\textit{(\alph*)}]
  \item automatic online PID adaptation and learning capability.
  \item adaptation only requires error and control signals.
  \item fast instant reaction to error changes at every iteration.
\end{enumerate*}

These characteristics make the NN-PID more robust and better at rejecting disturbances. On the other hand, the GA-PID produces optimized PID gains in terms of response speed, precision and overshoot but does not generalize well for different cases of disturbances and set-points. The NN-PID has been found to be more efficient in terms of adaptability and generalizes better. However, its performance significantly depends on weights initialization and choice of learning rate. In general, the GA-PID was smoother and faster than the NN-PID, while the latter was more adaptive and more robust. Future work should evaluate other optimization algorithms and address the smooth transition between different PID gains as well as explore other controllers. 

\bibliographystyle{ieeeconf}
\bibliography{library}

\begin{thebibliography}{99}

\bibitem{1} K. Osman, M. F. Rahmat, and M. A. Ahmad, “Modelling and controller design for a cruise control system,” Proc. 2009 5th Int. Colloq. Signal Process. Its Appl. CSPA 2009, no. 1, pp. 254–258, 2009


\bibitem{2} A. Simorgh, A. Marashian, and A. Razminia, “Adaptive PID Control Design for Longitudinal Velocity Control of Autonomous Vehicles,” Proc. - 2019 6th Int. Conf. Control. Instrum. Autom. ICCIA 2019, pp. 1–6, 2019

\bibitem{3} J. E. A. Dias, G. A. S. Pereira, and R. M. Palhares, “Longitudinal Model Identification and Velocity Control of an Autonomous Car,” IEEE Trans. Intell. Transp. Syst., vol. 16, no. 2, pp. 776–786, 2015

\bibitem{4} H. Kim, D. Kim, I. Shu, and K. Yi, “Time-varying parameter adaptive vehicle speed control,” IEEE Trans. Veh. Technol., vol. 65, no. 2, pp. 581–588, 2016

\bibitem{5} F. Gao, X. Hu, S. E. Li, K. Li, and Q. Sun, “Distributed Adaptive Sliding Mode Control of Vehicular Platoon with Uncertain Interaction Topology,” IEEE Trans. Ind. Electron., vol. 65, no. 8, pp. 6352–6361, 2018

\bibitem{6} S. Xu, H. Peng, Z. Song, K. Chen, and Y. Tang, “Accurate and Smooth Speed Control for an Autonomous Vehicle,” IEEE Intell. Veh. Symp. Proc., vol. 2018-June, no. Iv, pp. 1976–1982, 2018


\bibitem{7} F. Walz and S. Hohmann, “Model predictive longitudinal motion control for low velocities on known road profiles,” Veh. Syst. Dyn., vol. 58, no. 8, pp. 1310–1328, 2020



\bibitem{8} W. Farag and Z. Saleh, “Tuning of PID track followers for autonomous driving,” 2018 Int. Conf. Innov. Intell. Informatics, Comput. Technol. 3ICT 2018, pp. 1–7, 2018

\bibitem{9} A. A. M. Zahir, S. S. N. Alhady, W. A. F. W. Othman, and M. F. Ahmad, “Genetic algorithm optimization of pid controller for brushed DC motor,” Lect. Notes Mech. Eng., vol. 0, no. 9789811087875, pp. 427–437, 2018

\bibitem{10} M. K. Debnath, R. Agrawal, S. R. Tripathy, and S. Choudhury, “Artificial neural network tuned PID controller for LFC investigation including distributed generation,” Int. J. Numer. Model. Electron. Networks, Devices Fields, vol. 33, no. 5, pp. 1–17, 2020

\bibitem{11} G. Han, W. Fu, W. Wang, and Z. Wu, “The lateral tracking control for the intelligent vehicle based on adaptive PID neural network,” Sensors (Switzerland), vol. 17, no. 6, pp. 1–15, 2017

\bibitem{12} L. Nie, J. Guan, C. Lu, H. Zheng, and Z. Yin, “Longitudinal speed control of autonomous vehicle based on a self-adaptive PID of radial basis function neural network,” IET Intell. Transp. Syst., vol. 12, no. 6, pp. 485–494, 2018

\bibitem{13} B. Zhao, H. Wang, Q. Li, J. Li, and Y. Zhao. PID Trajectory Tracking Control of Autonomous Ground Vehicle Based on Genetic Algorithm. Proceedings of the 31st Chinese Control and Decision Conference, CCDC 2019, 3677–3682. 


\bibitem{14} Z. Sun and G. G. Zhu, Design and control of automotive propulsion systems, 1st Editio. Boca Raton: CRC Press, 2014.

\bibitem{15} A. Desreveaux, M. Ruba, A. Bouscayrol, G. M. Sirbu and C. Martis, "Comparisons of Models of Electric Drives for Electric Vehicles," 2019 IEEE Vehicle Power and Propulsion Conference (VPPC), Hanoi, Vietnam, 2019, pp. 1-5

\bibitem {16} Pacejka, H. B. (2008). Vehicle System Dynamics : International Journal of Vehicle Mechanics and Mobility. International Journal of Vehicle Mechanics and Mobility, (August 2012), 37–41.

\bibitem{17} A. Q. Ansari, Ibraheem, and S. Katiyar, “Application of ant colony algorithm for calculation and analysis of performance indices for adaptive control system,” Proc. Int. Conf. Innov. Appl. Comput. Intell. Power, Energy Control. with Their Impact Humanit. CIPECH 2014, no. November, pp. 466–471, 2014

\end{thebibliography}

\end{document}